\documentclass[12pt]{article}
\usepackage{amssymb,amsmath,amscd, amsthm,stmaryrd,verbatim, mathrsfs}
\numberwithin{equation}{section}

\begin{document}
\baselineskip=14pt

\newcommand{\la}{\langle}
\newcommand{\ra}{\rangle}
\newcommand{\psp}{\vspace{0.4cm}}
\newcommand{\pse}{\vspace{0.2cm}}
\newcommand{\ptl}{\partial}
\newcommand{\dlt}{\delta}
\newcommand{\sgm}{\sigma}
\newcommand{\al}{\alpha}
\newcommand{\be}{\beta}
\newcommand{\G}{\Gamma}
\newcommand{\gm}{\gamma}
\newcommand{\vs}{\varsigma}
\newcommand{\Lmd}{\Lambda}
\newcommand{\lmd}{\lambda}
\newcommand{\td}{\tilde}
\newcommand{\vf}{\varphi}
\newcommand{\yt}{Y^{\nu}}
\newcommand{\wt}{\mbox{wt}\:}
\newcommand{\rd}{\mbox{Res}}
\newcommand{\ad}{\mbox{ad}}
\newcommand{\stl}{\stackrel}
\newcommand{\ol}{\overline}
\newcommand{\ul}{\underline}
\newcommand{\es}{\epsilon}
\newcommand{\dmd}{\diamond}
\newcommand{\clt}{\clubsuit}
\newcommand{\vt}{\vartheta}
\newcommand{\ves}{\varepsilon}
\newcommand{\dg}{\dagger}
\newcommand{\tr}{\mbox{Tr}}
\newcommand{\ga}{{\cal G}({\cal A})}
\newcommand{\hga}{\hat{\cal G}({\cal A})}
\newcommand{\Edo}{\mbox{End}\:}
\newcommand{\for}{\mbox{for}}
\newcommand{\kn}{\mbox{ker}}
\newcommand{\Dlt}{\Delta}
\newcommand{\rad}{\mbox{Rad}}
\newcommand{\rta}{\rightarrow}
\newcommand{\mbb}{\mathbb}
\newcommand{\lra}{\Longrightarrow}
\newcommand{\X}{{\cal X}}
\newcommand{\Y}{{\cal Y}}
\newcommand{\Z}{{\cal Z}}
\newcommand{\U}{{\cal U}}
\newcommand{\V}{{\cal V}}
\newcommand{\W}{{\cal W}}
\newcommand{\sta}{\theta}
\setlength{\unitlength}{3pt}
\newcommand{\msr}{\mathscr}
\newcommand{\wht}{\widehat}
\newcommand{\mfk}{\mathfrak}

\begin{center}{\large \bf Congruence Classes of Supporting the Erd\"{o}s-Straus\\\pse Conjecture II: Wild Solutions} 
\end{center}

\vspace{0.2cm}

\begin{center}{\large  Xiaoping Xu}\end{center}

\pse

\begin{center}{
 HLM, Institute of Mathematics, Academy of Mathematics \& System
Sciences\\ Chinese Academy of Sciences, Beijing 100190, P.R. China
\\ \& School of Mathematics, University of Chinese Academy of Sciences,\\ Beijing 100049, P.R. China}\end{center}

\begin {abstract}
\quad

In 1948, Erd\"{o}s and Straus formulated a conjecture : for any positive integer $n>2$, there exist positive integers $n_1,n_2$ and $n_3$ such that
\begin{equation}\frac{4}{n}=\frac{1}{n_1}+\frac{1}{n_2}+\frac{1}{n_3},\nonumber\end{equation}
which is still open. It is known that the conjecture holds  if one can prove it for any prime $n\equiv 1\;(\mbox{mod}\;24)$. If $n=24m+1$ and $n_1\leq n_2,n_3$, then $n_1=6m+k$ with $1\leq k\leq 12m$. A solution $(n_1,n_2,n_3)$ of the above equation is called a {\it tame solution} if $n_2$ and $n_3$ are factors of $(6m+k)(24m+1)$. We call $n=24m+1$ {\it wild} if it does not have any tame solution. Based on the information in our earlier
 work on tame solutions posed in arXiv, Howerton found that there are only fourteen wild primes of the form $n=24m+1\leq 2.4\times 10^{11}$.
 In this paper, we derive thirty-four families of wild solutions of the above equation, which contain the solvability of the fourteen wild primes. Together with our earlier tame polynomial solutions, numeric test shows that they cover all the primes of the form $24m+1$.

 \vspace{0.3cm}

\noindent{\it Keywords}:\hspace{0.3cm} Erd\"{o}s-Straus conjecture; Egyptian fraction; congruence class; tame solution; wild solution; wild prime.
\end{abstract}

\section {Introduction}
Ancient Egyptians used sums of unit fractions (whose numerators are 1) to express fractions due to their ways of distributing food.
For instance,  $5/8$ was interpreted by ancient Egyptians as distributing five pancakes  fairly among eight people. They cut first four
pancakes into halves and then cut the last one into eight equal pieces. So each person got the same share: $1/2+1/8$ pancakes. This amazingly interpreted the mathematical equation
\begin{equation}\frac{5}{8}=\frac{1}{2}+\frac{1}{8}.\end{equation}
So an Egyptian fraction is a sum of distinct unit fractions. By repeatedly applying the simple fact
\begin{equation}\frac{1}{k}=\frac{1}{k+1}+\frac{1}{k(k+1)},\end{equation}
one can easily prove that any fraction is a Egyptian fraction. However,  it is very difficult to determine if a fraction can be expressed as a sum of fixed number of unit fractions (cf. \cite{ETao} for an excellent exposition and extensive references). For example, it is difficult to know if a fraction can be written as a sum of two unit fractions (e.g., cf. \cite{HV1,HV2,HV3,Jc2, Ps}). In 1948, Erd\"{o}s and Straus formulated a conjecture : for any positive integer $n>2$, there exist positive integers $n_1,n_2$ and $n_3$ such that
\begin{equation}\frac{4}{n}=\frac{1}{n_1}+\frac{1}{n_2}+\frac{1}{n_3},\end{equation}
which is still open up to now. For each $n$, the number of solution can go to infinity as $n$ does. Elscholtz and Tao \cite{ETao} found excellent bounds for it.

Mordell \cite{Ml} proved that the conjecture holds for positive integers in  834 congruence classes modulo 840. Terzi \cite{Tdj} used computer to verify that the conjecture holds for positive integers in  the congruence classes modulo 120120 except 198 classes. Kotsireas \cite{Ki}  verified the conjecture for every $n<10^{10}$. There are other interesting partial results on the conjecture or related works (e.g., cf. \cite{Bk, BL, Ht, IW, Ld, Ms, Ps, Rl, Sjw1, Sjw2, Sjw3, Sj, ST, Tdj, Vr, Ww}).

It is known that one only needs to prove the conjecture for any prime number $n$ such that $n\equiv 1\;(\mbox{mod}\;24)$. Suppose that $n=24m+1$ and $n_1\leq n_2,n_3$. Then $n_1=6m+k$ with $1\leq k\leq 12m$ and
\begin{equation}\frac{4}{24m+1}=\frac{1}{6m+k}+\frac{4k-1}{(6m+k)(24m+1)}.\end{equation}
If
\begin{equation}\frac{1}{n_2}=\frac{\Im_1}{(6m+k)(24m+1)}\quad\mbox{and}\quad \frac{1}{n_2}=\frac{\Im_2}{(6m+k)(24m+1)}\end{equation}for some positive integers $\Im_1$ and $\Im_2$ such that $\Im_1+\Im_2=4k-1$, we call the triple $(n_1,n_2,n_3)$ a {\it tame solution} of the Erd\"{o}s-Straus equation (1.3) for the positive integer $n=24m+1$. Moreover, we call $n=24m+1$ {\it wild} if it does not have any tame solution. Based on the information tion in our earlier
 work on tame solutions posed in arXiv, Howerton found that there are only fourteen wild primes of the form $n=24m+1\leq 2.4\times 10^{11}$.

By analysing the  obtained wild solutions, we find a way of obtaining them as follows. In (1.4), we multiply a proper integer $\ell \geq 2$ both to the numerator and denominator of the fraction
\begin{equation}\frac{4k-1}{(6m+k)(24m+1)}\lra\frac{\ell(4k-1)}{\ell(6m+k)(24m+1)}\end{equation} and find positive integers $\Im'_1$ and $\Im'_2$ such that
\begin{equation}\Im'_1+\Im'_2=\ell(4k-1),\end{equation}
\begin{equation} g.c.d(\ell,\Im'_1,\Im'_2)=1\end{equation} and
\begin{equation}\frac{1}{n_2}=\frac{\Im'_1}{\ell(6m+k)(24m+1)}\quad\mbox{and}\quad \frac{1}{n_2}=\frac{\Im_2'}{\ell(6m+k)(24m+1)}\end{equation}
are unit fractions. We call the integer $\ell$ a {\it multiplier} of the wild solution $(n_1,n_2,n_3)$. Moreover, we call $\Im'_1$ and $\Im'_2$
the {\it numerator summands} for $(n_1,n_2,n_3)$. According to our definitions, a tame prime of the form $24m+1$ may have wild solutions, but a wild prime of the form $24m+1$ only has wild solutions. Based on the information in our earlier
 work on tame solutions posed in arXiv, Howerton \cite{Hb} found that there are only fourteen wild primes of the form $n=24m+1\leq 2.4\times 10^{11}$. In this paper, we derive thirty-four families of wild solutions of the above equation according to (1,6)-(1.9), which contain the solvability of the fourteen wild primes. Together with our  tame polynomial solutions in \cite{Xx1}, numeric test shows that they cover all the primes of the form $24m+1$. The paper is organized as follows.\pse

In Section 2, we find a relatively simple wild solution for each of the fourteen wild primes of the form $24m+1$. In Section 3, we solve the Erd\"{o}s-Straus equation for $k\leq 5$ by the method described in (1.6)-(1.9). Moreover, we solve the equation for $k=6$ in Section 4.

\section{Solutions for the Fourteen Wild Primes}

Howerton \cite{Hb} found that there are only wild primes of the form $24m+1$ with $m\leq 10^9$. They are
\begin{equation}409,\quad 577,\quad 5569,\quad 9601,\quad 23929,\quad 83449,\quad 102001,\quad 329617,\quad 712321,\end{equation}
 \begin{equation}1134241,\quad 1724209,\quad 1726201,\quad 5212561,\quad 8813281.\end{equation}
In the form of $n=24m+1$, the corresponding $m$'s are as follows:
\begin{equation}17,\quad 24,\quad 232,\quad 400,\quad 997,\quad 3477,\quad 4250,\quad 13734,\quad
29680,\end{equation}
\begin{equation}47260,\quad 71842,\quad 71925,\quad217190,\quad 367220.\end{equation}
In this section, we derive a relatively simple wild solution for each of them as the starting point of this paper.

 Throughout this paper, we denote by $\mbb N$ the set of nonnegative integers. Suppose $m\equiv 4\;(\mbox{mod}\;5)$; that is, $m=4+5c$ for some $c\in\mbb N$. We have
\begin{eqnarray}\frac{4}{24m+1}&=&\frac{1}{6m+1}+\frac{3}{(6m+1)(24m+1)}\nonumber\\ &=& \frac{1}{6m+1}+\frac{6}{2(30c+25)(24m+1)}\nonumber\\ &=&
\frac{1}{6m+1}+\frac{5+1}{10(6c+5)(24m+1)}\nonumber\\ &=&\frac{1}{6m+1}+\frac{1}{2(6c+5)(24m+1)}+\frac{1}{10(6c+5)(24m+1)}.\end{eqnarray}
Taking $m=24$ and $c=4$, we get $24m+1=577$. $6m+1=145$ and
\begin{equation}\frac{4}{577}=\frac{1}{145}+\frac{1}{58\times 577}+\frac{1}{290\times 577}.\end{equation}
When $m=13734$, $24m+1=329617$. $6m+1=82405$, $c=2746$ and
\begin{equation}\frac{4}{329617}=\frac{1}{82405}+\frac{1}{32962\times 329617}+\frac{1}{164810\times 329617}.\end{equation}

Next we assume $m\equiv 3\;(\mbox{mod}\;7)$; that is, $m=3+7c$ for some $c\in\mbb N$. We have
\begin{eqnarray}\frac{4}{24m+1}&=&\frac{1}{6m+3}+\frac{11}{(6m+3)(24m+1)}\nonumber\\ &=& \frac{1}{6m+3}+\frac{22}{2(42c+21)(24m+1)}\nonumber\\ &=&
\frac{1}{6m+3}+\frac{21+1}{42(2c+1)(24m+1)}\nonumber\\ &=&\frac{1}{6m+3}+\frac{1}{2(2c+1)(24m+1)}+\frac{1}{42(2c+1)(24m+1)}.\end{eqnarray}
Taking $m=17$ and $c=2$, we get  $24m+1=409$. $6m+3=105$ and
\begin{equation}\frac{4}{409}=\frac{1}{105}+\frac{1}{10\times 409}+\frac{1}{210\times 409}.\end{equation}
When $m=47260$ and $c=6751$, we have $24m+1=1134241$, $6m+3=283563$ and
\begin{equation}\frac{4}{1134241}=\frac{1}{283563}+\frac{1}{27006\times 1134241}+\frac{1}{567126\times 1134241}.\end{equation}
If $m=367220$ and $c=52031$, we have $24m+1=8813281$, $6m+3=2203323$ and
\begin{equation}\frac{4}{8813281}=\frac{1}{2203323}+\frac{1}{104063\times 8813281}+\frac{1}{2185323\times 8813281}.\end{equation}

Suppose $m\equiv 11\;(\mbox{mod}\;17)$; that is, $m=11+17c$ for some $c\in\mbb N$. We have
\begin{eqnarray}\frac{4}{24m+1}&=&\frac{1}{6m+2}+\frac{7}{(6m+3)(24m+1)}\nonumber\\ &=& \frac{1}{6m+2}+\frac{35}{5(102c+68)(24m+1)}\nonumber\\ &=&
\frac{1}{6m+2}+\frac{34+1}{5\times 34(3c+2)(24m+1)}\nonumber\\ &=&\frac{1}{6m+2}+\frac{1}{5(3c+2)(24m+1)}+\frac{1}{170(3c+2)(24m+1)}.\end{eqnarray}
Taking $m=232$ and $c=13$, we get  $24m+1=5569$, $6m+2=1394$ and
\begin{equation}\frac{4}{5569}=\frac{1}{1394}+\frac{1}{1141545}+\frac{1}{38815930}
.\end{equation}
When $m=997$, we have $c=58$,  $24m+1=23929$, $6m+2=5984$ and
\begin{equation}\frac{4}{23929}=\frac{1}{5984}+\frac{1}{21057520} +\frac{1}{715955680}.\end{equation}

Let $m=400$. Then $24m+1=9601$ is a primes. Moreover, $6m+5=2405=5\times 13\times 37$. So
\begin{eqnarray}\frac{4}{9601}&=&\frac{1}{2405}+\frac{19}{2405\times 9601}\nonumber\\ &=&\frac{1}{2405}+\frac{19}{(5\times 13\times 37)\times 9601}\nonumber\\ &=&\frac{1}{2405}+\frac{38}{2\times(5\times 13\times 37)\times 9601}\nonumber\\ &=&\frac{1}{2405}+\frac{37+1}{2\times(5\times 13\times 37)\times 9601}\nonumber\\ &=&\frac{1}{2405}+\frac{1}{2\times(5\times 13)\times 9601}+\frac{1}{2\times(5\times 13\times 37)\times 9601}\nonumber\\ &=&\frac{1}{2405}+\frac{1}{1248130}+\frac{1}{46180810}
.\end{eqnarray}

In the case $m=3477$, $24m+1=83449$ is a primes. Moreover, $6m+3=20865=3\times 5\times 13\times 107$. Therefore,
\begin{eqnarray}\frac{4}{83449}&=&\frac{1}{20865}+\frac{11}{20865\times 83449}\nonumber\\ &=&\frac{1}{20865}+\frac{11}{(3\times 5\times 13\times 107)\times 83449}\nonumber\\ &=&\frac{1}{20865}+\frac{66}{6\times (3\times 65\times 107)\times 83449}\nonumber\\ &=&\frac{1}{20865}+\frac{65+1}{6\times (3\times 65\times 107)\times 83449}\nonumber\\ &=&\frac{1}{20865}+\frac{1}{6\times (3\times  107)\times 83449}+\frac{1}{6\times (3\times 65\times 107)\times 83449}\nonumber\\ &=&\frac{1}{20865}+\frac{1}{160722774}+\frac{1}{10446980310}.\end{eqnarray}

Assume $m=4250$, $24m+1=102001$ is a primes. Moreover, $6m+2=25502=2\times 41\times 311$. Thus
\begin{eqnarray}\frac{4}{102001}&=&\frac{1}{25502}+\frac{7}{25502\times 102001}\nonumber\\ &=&\frac{1}{25502}+\frac{7}{(2\times 41\times 311)\times 102001}\nonumber\\ &=&\frac{1}{25502}+\frac{42}{6\times(2\times 41\times 311)\times 102001}\nonumber\\ &=&\frac{1}{25502}+\frac{41+1}{6\times(2\times 41\times 311)\times 102001}\nonumber\\ &=&\frac{1}{25502}+\frac{1}{6\times(2\times 311)\times 102001}+\frac{1}{6\times(2\times 41\times 311)\times 102001}\nonumber\\ &=&\frac{1}{25502}+\frac{1}{380667732}+\frac{1}{15607377012}.\end{eqnarray}

We consider $m=29680$, $24m+1=712321$ is a primes. Moreover, $6m+6=178086=2\times 3\times 67\times 443$. Thus
\begin{eqnarray}\frac{4}{712321}&=&\frac{1}{178086}+\frac{23}{178086\times 712321}\nonumber\\ &=&\frac{1}{178086}+\frac{23}{(2\times 3\times 67\times 443)\times 712321}\nonumber\\ &=&\frac{1}{178086}+\frac{69}{3\times(2\times 3\times 67\times 443)\times 712321}\nonumber\\ &=&\frac{1}{178086}+\frac{67+2}{3\times(2\times 3\times 67\times 443)\times 712321}\nonumber\\ &=&\frac{1}{178086}+\frac{1}{3\times(2\times 3\times 443)\times 712321}+\frac{1}{3\times(3\times 67\times 443)\times 712321}\nonumber\\ &=&\frac{1}{178086}+\frac{1}{5680047654}+\frac{1}{190281596409}.\end{eqnarray}

Now we assume $m\equiv  9\;\;(\mbox{mod}\;29);$ equivalently, $m=9+29c$ for some $c\in\mbb N.$
In this case,
\begin{eqnarray}\frac{4}{24m+1}&=&\frac{1}{6m+4}+\frac{15}{2(3m+2)(24m+1)}\nonumber\\&=&
\frac{1}{6m+4}+\frac{30}{4(3m+2)(24m+1)}\nonumber\\&=&
\frac{1}{6m+4}+\frac{29+1}{116(3c+1)(24m+1)}\nonumber\\&=&
\frac{1}{6m+4}+\frac{1}{4(3c+1)(24m+1)}+\frac{1}{116(3c+1)(24m+1)}
.\end{eqnarray}
Taking $m=71842$, we get $c=2477$ and $24m+1=1724209$ is a wild prime. By the above equation,
\begin{eqnarray}\frac{4}{1724209}&=&
\frac{1}{431056}+\frac{1}{29728\times 1724209}+\frac{1}{862112\times 1724209}
.\end{eqnarray}
When $m=217190$, we have $c=7489$ and $24m+1=5212561$ is a wild prime. According to (2.18).,
\begin{eqnarray}\frac{4}{5212561}&=&\frac{1}{1303144}+\frac{1}{89872\times 5212561}+\frac{1}{2606288\times 5212561}
.\end{eqnarray}

Finally we consider $m=71925$, $24m+1=1726201$ is a wild primes. Moreover, $6m+6=431556=4\times 3\times 35963$. In fact,
\begin{equation}3\times 35963=107889,\quad 23\times 4691=107893=3\times 35963+4.\end{equation}
 Thus
\begin{eqnarray}\frac{4}{1726201}&=&\frac{1}{431556}+\frac{23}{431556\times 1726201}\nonumber\\ &=&
\frac{1}{431556}+\frac{4691\times 23}{4691\times 431556\times 1726201}\nonumber\\ &=&
\frac{1}{431556}+\frac{3\times 35963+4}{4691\times 4\times 3\times 35963\times 1726201}\nonumber\\ &=&
\frac{1}{431556}+\frac{1}{18764\times 1726201}
+\frac{1}{4691\times 107889\times 1726201}.\end{eqnarray}

\pse

The above examples suggest us to solve the Erd\"{o}s-Straus equation (1.3) for wild solutions under the conditions $k\leq 6$. In fact, we solve the equation by the method described in (1.6)-(1.9) and under the condition
\begin{equation}\Im'_1|[\ell(6m+k)],\quad \Im'_2|[\ell(6m+k)].\end{equation}
In the tame case, the above expression holds as long as we assume that $n=24m+1$ is a prime (cf. \cite{Xx1}). Computer shows that it does not hold in general wild case.

\section{Wild Solutions with $k\leq 5$}

In this section, we present the wild solutions for an integer $n=24m+1$ with $k=1,2,3,4,5$, where
\begin{equation}\frac{4}{n}=\frac{1}{6m+k}+\frac{4k-1}{(6m+k+1)(24m+1)}.\end{equation}Since the case $k=6$ is extremely complicated, we will solve it in next section.

\subsection{Case $k=1$}

In this case, we consider
\begin{eqnarray}\frac{4}{n}&=&\frac{1}{6m+1}+\frac{3}{(6m+1)(24m+1)}\nonumber\\&=&
\frac{1}{6m+1}+\frac{2\ell\times 3}{2\ell\times(6m+1)(24m+1)}\nonumber\\&=&\frac{1}{6m+1}+\frac{(6\ell-1)+1}{2\ell(6m+1)(24m+1)},\end{eqnarray}where $\ell$ is a positive integer.
Motivated from (2.5), we assume
\begin{equation}6m+1\equiv 0\quad(\mbox{mod}\;6\ell-1).\end{equation}
Equivalently,
\begin{equation}m+\ell\equiv 0\quad(\mbox{mod}\;6\ell-1).\end{equation}
Thus \begin{equation}m\equiv 5\ell-1\quad(\mbox{mod}\;6\ell-1).\end{equation}
When $\ell=1$, we recover the condition $m\equiv 4\;(\mbox{mod}\;5)$ for (2.5). Hence
\begin{equation}m=5\ell-1+c(6\ell-1)\qquad\mbox{for some}\;\;c\in\mbb N.\end{equation}

Suppose that $m$ is of the above form. Then
\begin{equation}6m+1=30\ell-6+6c(6\ell-1)+1=(6c+5)(6\ell-1).\end{equation}According to (3.2),
\begin{eqnarray}\frac{4}{24m+1}&=&\frac{1}{6m+1}+\frac{(6\ell-1)+1}{2\ell(6m+1)(24m+1)}\nonumber\\&=&
\frac{1}{6m+1}+\frac{(6\ell-1)+1}{2\ell(6c+5)(6\ell-1)(24m+1)}\nonumber\\&=&\frac{1}{6m+1}+\frac{1}{2\ell (6c+5)(24m+1)}+\frac{1}{2\ell (6c+5)(6\ell-1)(24m+1)}.\end{eqnarray}\pse

{\bf Theorem 3.1}\quad {\it For any positive integer $m$ of the form (3.6), we have the wild solution (3.8) of the Erd\"{o}s-Straus equation.}\psp

\subsection{Case $k=2$}

In this case, we first consider
\begin{eqnarray}\frac{4}{n}&=&\frac{1}{6m+2}+\frac{7}{(6m+2)(24m+1)}\nonumber\\&=&
\frac{1}{6m+2}+\frac{2\ell\times 7}{2\ell\times(6m+2)(24m+1)}\nonumber\\&=&\frac{1}{6m+2}+\frac{(14\ell-1)+1}{4\ell(3m+1)(24m+1)},\end{eqnarray}where $\ell$ is a positive integer.
Motivated from (2.5), we assume
\begin{equation}3m+1\equiv 0\quad(\mbox{mod}\;14\ell-1).\end{equation}

It is impossible if $\ell\equiv 2\;(\mbox{mod}\;3)$. Assume $\ell=3r$ with $0<r\in\mbb N$, Then the above equation becomes
\begin{equation}3m+1\equiv 0\quad(\mbox{mod}\;42r-1).\end{equation}
Equivalently,
\begin{equation}m+14r\equiv 0\quad(\mbox{mod}\;42r-1).\end{equation}
Thus \begin{equation}m\equiv 28r-1\quad(\mbox{mod}\;42r-1).\end{equation}
So
\begin{equation}m=28r-1+c(42r-1)\qquad\mbox{for some}\;\;c\in\mbb N.\end{equation}

Suppose that $m$ is of the above form. Then
\begin{equation}3m+1=84r-3+3c(42r-1)+1=(3c+2)(42r-1).\end{equation}According to (3.9),
\begin{eqnarray}&&\frac{4}{24m+1}=\frac{1}{6m+2}+\frac{(14\ell-1)+1}{4\ell(3m+1)(24m+1)}\nonumber\\&=&
\frac{1}{6m+2}+\frac{(42r-1)+1}{12r(3c+2)(42r-1)(24m+1)}\nonumber\\&=&
\frac{1}{6m+2}+\frac{1}{12r(3c+2)(24m+1)}+\frac{1}{12r(3c+2)(42r-1)(24m+1)}.\end{eqnarray}\pse
When $r=1$ and $c=103$, the above equation yields (2.17).

Next we assume $\ell=3r+1$ with $r\in\mbb N$, Then Equation (3.10) becomes
\begin{equation}3m+1\equiv 0\quad(\mbox{mod}\;42r+13).\end{equation}
Equivalently,
\begin{equation}m\equiv 14r+4\quad(\mbox{mod}\;42r+13).\end{equation}
Hence
\begin{equation}m=14r+4+c(42r+13)\qquad\mbox{for some}\;\;c\in\mbb N.\end{equation}

If $m$ is of the above form, then
\begin{equation}3m+1=42r+12+3c(42r+13)+1=(3c+1)(42r+13).\end{equation}By (3.9),
\begin{eqnarray}\frac{4}{24m+1}&=&\frac{1}{6m+2}+\frac{(14\ell-1)+1}{4\ell(3m+1)(24m+1)}\nonumber\\&=&
\frac{1}{6m+2}+\frac{(42r+13)+1}{4(3r+1)(3c+1)(42r+13)(24m+1)}\nonumber\\&=&
\frac{1}{6m+2}+\frac{1}{4(3r+1)(3c+1)(24m+1)}\nonumber\\&&+\frac{1}{4(3r+1)(3c+1)(42r+13)(24m+1)}.\end{eqnarray}\pse

Now we consider
\begin{eqnarray}\frac{4}{n}&=&\frac{1}{6m+2}+\frac{7}{(6m+2)(24m+1)}\nonumber\\&=&
\frac{1}{6m+2}+\frac{(2\ell+1)\times 7}{(2\ell+1)\times(6m+2)(24m+1)}\nonumber\\&=&\frac{1}{6m+2}+\frac{(14\ell+6)+1}{2(2\ell+1)(3m+1)(24m+1)},\end{eqnarray}where $\ell$ is a positive integer.
Motivated from (2.5), we assume
\begin{equation}3m+1\equiv 0\quad(\mbox{mod}\;7\ell+3).\end{equation}
It is impossible if $\ell\equiv 0\;(\mbox{mod}\;3)$. Assume $\ell=3r+1$ with $r\in\mbb N$, Then the above equation becomes
\begin{equation}3m+1\equiv 0\quad(\mbox{mod}\;21r+10).\end{equation}
Equivalently,
\begin{equation}m\equiv 7r+3\quad(\mbox{mod}\;21r+10).\end{equation}
Therefore,
\begin{equation}m=7r+3+c(21r+10)\qquad\mbox{for some}\;\;c\in\mbb N.\end{equation}

When $m$ is of the above form,
\begin{equation}3m+1=21r+9+3c(21r+10)+1=(3c+1)(21r+10).\end{equation}By (3.22),
\begin{eqnarray}\frac{4}{24m+1}&=&\frac{1}{6m+2}+\frac{(14\ell+6)+1}{2(2\ell+1)(3m+1)(24m+1)}
\nonumber\\&=&\frac{1}{6m+2}+\frac{2(21r+10)+1}{6(2r+1)(3c+1)(21r+10)(24m+1)}
\nonumber\\&=&\frac{1}{6m+2}+\frac{1}{3(2r+1)(3c+1)(24m+1)}\nonumber\\&&+\frac{1}{6(2r+1)(3c+1)(21r+10)(24m+1)}
.\end{eqnarray}\pse

Next we assume $\ell=3r+2$ with $r\in\mbb N$, Then Equation (3.23) becomes
\begin{equation}3m+1\equiv 0\quad(\mbox{mod}\;21r+17).\end{equation}
Equivalently,
\begin{equation}m+7r+6\equiv 0\quad(\mbox{mod}\;21r+17);\end{equation}that is
\begin{equation}m\equiv 14r+11\quad(\mbox{mod}\;21r+17);\end{equation}
Hence
\begin{equation}m=14r+11+c(21r+17)\qquad\mbox{for some}\;\;c\in\mbb N.\end{equation}

Suppose that $m$ is of the above form. Then
\begin{equation}3m+1=42r+33+3c(21r+17)+1=(3c+2)(21r+17).\end{equation}By (3.22),
\begin{eqnarray}\frac{4}{24m+1}&=&\frac{1}{6m+2}+\frac{(14\ell+6)+1}{2(2\ell+1)(3m+1)(24m+1)}
\nonumber\\&=&\frac{1}{6m+2}+\frac{2(21r+17)+1}{2(6r+5)(3c+2)(21r+17)(24m+1)}
\nonumber\\&=&\frac{1}{6m+2}+\frac{1}{(6r+5)(3c+2)(24m+1)}\nonumber\\&&+\frac{1}{2(6r+5)(3c+2)(21r+17)(24m+1)}.\end{eqnarray}
When $r=0$, the above equation becomes (2.12).

Finally we consider
\begin{eqnarray}\frac{4}{n}&=&\frac{1}{6m+2}+\frac{7}{(6m+2)(24m+1)}\nonumber\\&=&
\frac{1}{6m+2}+\frac{(2\ell+1)\times 7}{(2\ell+1)\times(6m+2)(24m+1)}\nonumber\\&=&\frac{1}{6m+2}+\frac{(14\ell+5)+2}{2(2\ell+1)(3m+1)(24m+1)},\end{eqnarray}where $\ell$ is a positive integer.
Motivated from (2.15), we assume
\begin{equation}3m+1\equiv 0\quad(\mbox{mod}\;14\ell+5).\end{equation}
It is impossible if $\ell\equiv 2\;(\mbox{mod}\;3)$. Assume $\ell=3r$ with $0<r\in\mbb N$, Then the above equation becomes
\begin{equation}3m+1\equiv 0\quad(\mbox{mod}\;42r+5).\end{equation}
Equivalently,
\begin{equation}m+14r+2\equiv 0\quad(\mbox{mod}\;42r+5);\end{equation}that is
\begin{equation}m\equiv 28r+3\quad(\mbox{mod}\;42r+5);\end{equation}
Thus
\begin{equation}m=28r+3+c(42r+5)\qquad\mbox{for some}\;\;c\in\mbb N.\end{equation}

If $m$ is of the above form, then
\begin{equation}3m+1=84r+9+3c(42r+5)+1=(3c+2)(42r+5).\end{equation}By (3.35),
\begin{eqnarray}\frac{4}{24m+1}&=&\frac{1}{6m+2}+\frac{(14\ell+5)+2}{2(2\ell+1)(3m+1)(24m+1)}
\nonumber\\&=&\frac{1}{6m+2}+\frac{(42r+5)+2}{2(6r+1)(3c+2)(42r+5)(24m+1)}
\nonumber\\&=&\frac{1}{6m+2}+\frac{1}{2(6r+1)(3c+2)(24m+1)}\nonumber\\&&+\frac{1}{(6r+1)(3c+2)(42r+5)(24m+1)}.\end{eqnarray}\pse

Next we assume $\ell=3r+1$ with $r\in\mbb N$, Then Equation (3.36) becomes
\begin{equation}3m+1\equiv 0\quad(\mbox{mod}\;42r+19).\end{equation}
Equivalently,
\begin{equation}m\equiv 14r+6\quad(\mbox{mod}\;42r+19).\end{equation}
In this way,
\begin{equation}m=14r+6+c(42r+19)\qquad\mbox{for some}\;\;c\in\mbb N.\end{equation}

When $m$ is of the above form,
\begin{equation}3m+1=42r+18+3c(42r+19)+1=(3c+1)(42r+19).\end{equation}According to (3.35),
\begin{eqnarray}\frac{4}{24m+1}&=&\frac{1}{6m+2}+\frac{(14\ell+5)+2}{2(2\ell+1)(3m+1)(24m+1)}
\nonumber\\&=&\frac{1}{6m+2}+\frac{(42r+19)+2}{6(2r+1)(3c+1)(42r+19)(24m+1)}
\nonumber\\&=&\frac{1}{6m+2}+\frac{1}{6(2r+1)(3c+1)(24m+1)}\nonumber\\&&+\frac{1}{3(2r+1)(3c+1)(42r+19)(24m+1)}.\end{eqnarray}\pse

{\bf Theorem 3.2}\quad {\it For any positive integer $m$ of the form (3.14), we have the wild solution (3.16) of the Erd\"{o}s-Straus equation.
If $m$ is of the form (3.19), we have the wild solution (3.21) of the Erd\"{o}s-Straus equation. When $m$ is of the form (3.26), we have the wild solution (3.28) of the Erd\"{o}s-Straus equation. Suppose that $m$ is of the form (3.32), we have the wild solution (3.34) of the Erd\"{o}s-Straus equation. Let $m$ be of the form (3.40), we have the wild solution (3.42) of the Erd\"{o}s-Straus equation. Assume that
$m$ is of the form (3.45), we have the wild solution (3.47) of the Erd\"{o}s-Straus equation.
}\psp

\subsection{Case $k=3$}

In this case, we first consider
\begin{eqnarray}\frac{4}{n}&=&\frac{1}{6m+3}+\frac{11}{(6m+3)(24m+1)}\nonumber\\&=&
\frac{1}{6m+3}+\frac{2\ell\times 11}{6\ell\times(2m+1)(24m+1)}\nonumber\\&=&\frac{1}{6m+3}+\frac{(22\ell-1)+1}{6\ell(2m+1)(24m+1)},\end{eqnarray}where $\ell$ is a positive integer.
Suppose $\ell\equiv 1\;(\mbox{mod}\;3)$; that is, $\ell=3r+1$ for some $r\in\mbb N$. Then
\begin{equation}22\ell-1=66r+21=3(22r+7).\end{equation}
By (2.16), we assume
\begin{equation}2m+1\equiv 0\quad(\mbox{mod}\;22r+7).\end{equation}
Equivalently,
\begin{equation}m\equiv 11r+3\quad(\mbox{mod}\;22r+7).\end{equation}
Thus
\begin{equation}m=11r+3+c(22r+7)\qquad\mbox{for some}\;\;c\in\mbb N.\end{equation}

Assume that $m$ is of the above form. Then
\begin{equation}2m+1=22r+6+2c(22r+7)+1=(2c+1)(22r+7).\end{equation}
According to (3.48),
\begin{eqnarray}&&\frac{4}{24m+1}=\frac{1}{6m+3}+\frac{(22\ell-1)+1}{6\ell(2m+1)(24m+1)}
\nonumber\\&=&\frac{1}{6m+3}+\frac{3(22r+7)+1}{6(3r+1)(2c+1)(22r+7)(24m+1)}
\nonumber\\&=&\frac{1}{6m+3}+\frac{1}{2(3r+1)(2c+1)(24m+1)}\nonumber\\&&+\frac{1}{6(3r+1)(2c+1)(22r+7)(24m+1)}.\end{eqnarray}
When $r=0$, the above equation yields (2.8). A subtle point is that we have well used the factor ``3" in $6m+3$ of the denominator. \psp

Next we assume $\ell\equiv 0\;(\mbox{mod}\;3)$; that is, $\ell=3r$ for some $0<r\in\mbb N$. Then
\begin{equation}22\ell-1=66r-1.\end{equation}
Motivated from (2.103, we impose
\begin{equation}2m+1\equiv 0\quad(\mbox{mod}\;66r-1).\end{equation}
Equivalently,
\begin{equation}m\equiv 33r-1\quad(\mbox{mod}\;66r-1).\end{equation}
Thus
\begin{equation}m=33r-1+c(66r-1)\qquad\mbox{for some}\;\;c\in\mbb N.\end{equation}

When $m$ is of the above form,
\begin{equation}2m+1=66r-2+2c(66r-1)+1=(2c+1)(66r-1).\end{equation}
By (3.48),
\begin{eqnarray}&&\frac{4}{24m+1}=\frac{1}{6m+3}+\frac{(22\ell-1)+1}{6\ell(2m+1)(24m+1)}
\nonumber\\&=&\frac{1}{6m+3}+\frac{(66r-1)+1}{18r(2c+1)(66r-1)(24m+1)}
\nonumber\\&=&\frac{1}{6m+3}+\frac{1}{18r(2c+1)(24m+1)}\nonumber\\&&+\frac{1}{18r(2c+1)(66r-1)(24m+1)}.\end{eqnarray}
When $r=1$ and $c=53$, the above equation yields (2.13).\psp

Thirdly we assume $\ell\equiv 2\;(\mbox{mod}\;3)$; that is, $\ell=3r+2$ for some $r\in\mbb N$. Then
\begin{equation}22\ell-1=66r+43.\end{equation}
Motivated from (2.13), we impose
\begin{equation}2m+1\equiv 0\quad(\mbox{mod}\;66r+43).\end{equation}
Equivalently,
\begin{equation}m\equiv 33r+21\quad(\mbox{mod}\;66r+43).\end{equation}
Thus
\begin{equation}m=33r+21+c(66r+43)\qquad\mbox{for some}\;\;c\in\mbb N.\end{equation}

If $m$ is of the above form,
\begin{equation}2m+1=66r+42+2c(66r+43)+1=(2c+1)(66r+43).\end{equation}
Expression (3.48) yields
\begin{eqnarray}&&\frac{4}{24m+1}=\frac{1}{6m+3}+\frac{(22\ell-1)+1}{6\ell(2m+1)(24m+1)}
\nonumber\\&=&\frac{1}{6m+3}+\frac{(66r+43)+1}{6(3r+2)(2c+1)(66r+43)(24m+1)}
\nonumber\\&=&\frac{1}{6m+3}+\frac{1}{6(3r+2)(2c+1)(24m+1)}
\nonumber\\&&+\frac{1}{6(3r+2)(2c+1)(66r+43)(24m+1)}.\end{eqnarray}\pse

 Now we  consider
\begin{eqnarray}\frac{4}{n}&=&\frac{1}{6m+3}+\frac{11}{(6m+3)(24m+1)}\nonumber\\&=&
\frac{1}{6m+3}+\frac{2\ell\times 11}{6\ell\times(2m+1)(24m+1)}\nonumber\\&=&\frac{1}{6m+3}+\frac{(22\ell-3)+3}{6\ell(2m+1)(24m+1)},\end{eqnarray}where $\ell$ is a positive integer.
According to (1.8) with $\Im'_1=22\ell-3$ and $\Im'_2=3$, $\ell\not\equiv 0\;(\mbox{mod}\;3)$. Suppose $\ell\equiv 1\;(\mbox{mod}\;3)$; that is, $\ell=3r+1$ for some $0<r\in\mbb N$. Then
\begin{equation}22\ell-3=66r+19.\end{equation}
Motivated from (2.8), we assume
\begin{equation}2m+1\equiv 0\quad(\mbox{mod}\;66r+19).\end{equation}
Equivalently,
\begin{equation}m\equiv 33r+9\quad(\mbox{mod}\;66r+19).\end{equation}
Thus
\begin{equation}m=33r+9+c(66r+19)\qquad\mbox{for some}\;\;c\in\mbb N.\end{equation}

Assume that $m$ is of the above form. Then
\begin{equation}2m+1=66r+18+2c(66r+19)+1=(2c+1)(66r+19).\end{equation}
According to (3.67),
\begin{eqnarray}&&\frac{4}{24m+1}=\frac{1}{6m+3}+\frac{(22\ell-3)+3}{6\ell(2m+1)(24m+1)}
\nonumber\\&=&\frac{1}{6m+3}+\frac{(66r+19)+3}{6(3r+1)(2c+1)(66r+19)(24m+1)}
\nonumber\\&=&\frac{1}{6m+3}+\frac{1}{6(3r+1)(2c+1)(24m+1)}\nonumber\\&&+\frac{1}{2(3r+1)(2c+1)(66r+19)(24m+1)}.\end{eqnarray}\pse

Secondly we suppose $\ell\equiv 2\;(\mbox{mod}\;3)$; that is, $\ell=3r+2$ for some $0<r\in\mbb N$. Then
\begin{equation}22\ell-3=66r+41.\end{equation}
Motivated from (2.8), we assume
\begin{equation}2m+1\equiv 0\quad(\mbox{mod}\;66r+41).\end{equation}
Equivalently,
\begin{equation}m\equiv 33r+20\quad(\mbox{mod}\;66r+41).\end{equation}
Thus
\begin{equation}m=33r+20+c(66r+41)\qquad\mbox{for some}\;\;c\in\mbb N.\end{equation}

Let $m$ be of the above form. Then
\begin{equation}2m+1=66r+40+2c(66r+41)+1=(2c+1)(66r+41).\end{equation}
By (3.67),
\begin{eqnarray}&&\frac{4}{24m+1}=\frac{1}{6m+3}+\frac{(22\ell-3)+3}{6\ell(2m+1)(24m+1)}
\nonumber\\&=&\frac{1}{6m+3}+\frac{(66r+41)+3}{6(3r+1)(2c+1)(66r+41)(24m+1)}
\nonumber\\&=&\frac{1}{6m+3}+\frac{1}{6(3r+1)(2c+1)(24m+1)}\nonumber\\&&+\frac{1}{2(3r+1)(2c+1)(66r+41)(24m+1)}.\end{eqnarray}

\pse

{\bf Theorem 3.3}\quad {\it For any positive integer $m$ of the form (3.52), we have the wild solution (3.54) of the Erd\"{o}s-Straus equation.
If $m$ is of the form (3.58), we have the wild solution (3.60) of the Erd\"{o}s-Straus equation. When $m$ is of the form (3.64), we have the wild solution (3.66) of the Erd\"{o}s-Straus equation. Suppose that $m$ is of the form (3.71), we have the wild solution (3.73) of the Erd\"{o}s-Straus equation. Let $m$ be of the form (3.77), we have the wild solution (3.79) of the Erd\"{o}s-Straus equation.}\psp

\subsection{Case $k=4$}

In this case, we first consider
\begin{eqnarray}\frac{4}{n}&=&\frac{1}{6m+4}+\frac{15}{(6m+4)(24m+1)}\nonumber\\&=&
\frac{1}{6m+4}+\frac{2\ell\times 15}{2\ell\times2(3m+2)(24m+1)}\nonumber\\&=&\frac{1}{6m+4}+\frac{(30\ell-1)+1}{4\ell(3m+2)(24m+1)},\end{eqnarray}where $\ell$ is a positive integer.
Motivated from (2.8), we assume
\begin{equation}3m+2\equiv 0\quad(\mbox{mod}\;30\ell-1).\end{equation}
Equivalently,
\begin{equation}m\equiv 10\ell-1\quad(\mbox{mod}\;30\ell-1).\end{equation}
Thus
\begin{equation}m=10\ell-1+c(30\ell-1)\qquad\mbox{for some}\;\;c\in\mbb N.\end{equation}

Assume that $m$ is of the above form. Then
\begin{equation}3m+2=30\ell-3+3c(30\ell-1)+2=(3c+1)(30\ell-1).\end{equation}
According to (3.80),
\begin{eqnarray}\frac{4}{24m+1}&=&\frac{1}{6m+4}+\frac{(30\ell-1)+1}{4\ell(3m+2)(24m+1)}\nonumber\\&=&
\frac{1}{6m+4}+\frac{(30\ell-1)+1}{4\ell(3c+1)(30\ell-1)(24m+1)}\nonumber\\&=&
\frac{1}{6m+4}+\frac{1}{4\ell(3c+1)(24m+1)}\nonumber\\&&+\frac{1}{4\ell(3c+1)(30\ell-1)(24m+1)}
.\end{eqnarray}
Taking $\ell=1$, we get (2.19).

\pse

Next we  consider
\begin{eqnarray}\frac{4}{n}&=&\frac{1}{6m+4}+\frac{15}{(6m+4)(24m+1)}\nonumber\\&=&
\frac{1}{6m+4}+\frac{(2\ell+1)\times 15}{(2\ell+1)\times2(3m+2)(24m+1)}\nonumber\\&=&\frac{1}{6m+4}+\frac{(30\ell+14)+1}{2(2\ell+1)(3m+2)(24m+1)},\end{eqnarray}where $\ell$ is a positive integer.
Motivated from (2.8), we assume
\begin{equation}3m+2\equiv 0\quad(\mbox{mod}\;15\ell+7).\end{equation}
Equivalently,
\begin{equation}m\equiv 10\ell+4\quad(\mbox{mod}\;15\ell+7).\end{equation}
Thus
\begin{equation}m=10\ell+4+c(15\ell+7)\qquad\mbox{for some}\;\;c\in\mbb N.\end{equation}

When $m$ is of the above form,
\begin{equation}3m+2=30\ell+12+3c(15\ell+7)+2=(3c+2)(15\ell+7).\end{equation}
By (3.86),
\begin{eqnarray}\frac{4}{24m+1}&=&\frac{1}{6m+4}+\frac{(30\ell+14)+1}{2(2\ell+1)(3m+2)(24m+1)}
\nonumber\\&=&\frac{1}{6m+4}+\frac{2(15\ell+7)+1}{2(2\ell+1)(3c+2)(15\ell+7)(24m+1)}
\nonumber\\&=&\frac{1}{6m+4}+\frac{1}{(2\ell+1)(3c+2)(24m+1)}\nonumber\\&&+\frac{1}{2(2\ell+1)(3c+2)(15\ell+7)(24m+1)}.\end{eqnarray}
\pse

Finally we  consider
\begin{eqnarray}\frac{4}{n}&=&\frac{1}{6m+4}+\frac{15}{(6m+4)(24m+1)}\nonumber\\&=&
\frac{1}{6m+4}+\frac{(2\ell+1)\times 15}{(2\ell+1)\times2(3m+2)(24m+1)}\nonumber\\&=&\frac{1}{6m+4}+\frac{(30\ell+13)+2}{2(2\ell+1)(3m+2)(24m+1)},\end{eqnarray}where $\ell$ is a positive integer.
Motivated from (2.15), we assume
\begin{equation}3m+2\equiv 0\quad(\mbox{mod}\;30\ell+13).\end{equation}
Equivalently,
\begin{equation}m\equiv 20\ell+8\quad(\mbox{mod}\;30\ell+13).\end{equation}
Thus
\begin{equation}m=20\ell+8+c(30\ell+13)\qquad\mbox{for some}\;\;c\in\mbb N.\end{equation}

If $m$ is of the above form,
\begin{equation}3m+2=60\ell+24+3c(30\ell+13)+2=(3c+2)(30\ell+13).\end{equation}
Expression (3.92) gives
\begin{eqnarray}\frac{4}{24m+1}&=&\frac{1}{6m+4}+\frac{(30\ell+13)+2}{2(2\ell+1)(3m+2)(24m+1)}
\nonumber\\&=&\frac{1}{6m+4}+\frac{(30\ell+13)+2}{2(2\ell+1)(3c+2)(30\ell+13)(24m+1)}
\nonumber\\&=&\frac{1}{6m+4}+\frac{1}{2(2\ell+1)(3c+2)(24m+1)}\nonumber\\&&+\frac{1}{(2\ell+1)(3c+2)(30\ell+13)(24m+1)}.\end{eqnarray}
\pse

{\bf Theorem 3.4}\quad {\it For any positive integer $m$ of the form (3.83), we have the wild solution (3.85) of the Erd\"{o}s-Straus equation.
If $m$ is of the form (3.89), we have the wild solution (3.91) of the Erd\"{o}s-Straus equation. When $m$ is of the form (3.95), we have the wild solution (3.97) of the Erd\"{o}s-Straus equation. }\psp

\subsection{Case $k=5$}

In this case, we consider
\begin{eqnarray}\frac{4}{n}&=&\frac{1}{6m+5}+\frac{19}{(6m+5)(24m+1)}\nonumber\\&=&
\frac{1}{6m+5}+\frac{2\ell\times 19}{2\ell\times(6m+5)(24m+1)}\nonumber\\&=&\frac{1}{6m+5}+\frac{(38\ell-1)+1}{2\ell(6m+5)(24m+1)},\end{eqnarray}where $\ell$ is a positive integer.
As the earlier cases, we assume
\begin{equation}6m+5\equiv 0\quad(\mbox{mod}\;38\ell-1).\end{equation}
It is impossible if $\ell\equiv 2\;(\mbox{mod}\;3)$. Assume $\ell=3r$ with $0<r\in\mbb N$, Then the above equation becomes
\begin{equation}6m+5\equiv 0\quad(\mbox{mod}\;114r-1).\end{equation}
Equivalently,
\begin{equation}m+95r\equiv 0\quad(\mbox{mod}\;114r-1).\end{equation}
Thus \begin{equation}m\equiv 19r-1\quad(\mbox{mod}\;114r-1).\end{equation}
So
\begin{equation}m=19r-1+c(114r-1)\qquad\mbox{for some}\;\;c\in\mbb N.\end{equation}

When $m$ is of the above form,
\begin{equation}6m+5=114r-6+6c(114r-1)+5=(6c+1)(114r-1).\end{equation}
According to (3.98),
\begin{eqnarray}\frac{4}{24m+1}&=&\frac{1}{6m+5}+\frac{(38\ell-1)+1}{2\ell(6m+5)(24m+1)}\nonumber\\&=&
\frac{1}{6m+5}+\frac{(114r-1)+1}{6r(6c+1)(114r-1)(24m+1)}\nonumber\\&=&
\frac{1}{6m+5}+\frac{1}{6r(6c+1)(24m+1)}\nonumber\\&&+\frac{1}{6r(6c+1)(114r-1)(24m+1)}.\end{eqnarray}\pse

Suppose $\ell\equiv 1\;(\mbox{mod}\;3)$; that is, $\ell=3r+1$ for some $r\in\mbb N$. Then
\begin{equation}38\ell-1=114r+37.\end{equation}
Equation (3.99) becomes
\begin{equation}6m+5\equiv 0\quad(\mbox{mod}\;114r+37).\end{equation}
Equivalently,
 \begin{equation}m\equiv 95r+30\quad(\mbox{mod}\;114r+37).\end{equation}
So
\begin{equation}m=95r+30+c(114r+37)\qquad\mbox{for some}\;\;c\in\mbb N.\end{equation}

Let $m$ be of the above form,
\begin{equation}6m+5=570r+180+6c(114r+37)+5=(6c+5)(114r+37).\end{equation}
According to (3.98),
\begin{eqnarray}\frac{4}{24m+1}&=&\frac{1}{6m+5}+\frac{(38\ell-1)+1}{2\ell(6m+5)(24m+1)}\nonumber\\&=&
\frac{1}{6m+5}+\frac{(114r+37)+1}{2(3r+1)(6c+5)(114r+37)(24m+1)}\nonumber\\&=&
\frac{1}{6m+5}+\frac{1}{2(3r+1)(6c+5)(24m+1)}\nonumber\\&&+\frac{1}{2(3r+1)(6c+5)(114r+37)(24m+1)}.\end{eqnarray}
When $r=0$ and $c=10$, the above equation gives (2.15).\psp

{\bf Theorem 3.5}\quad {\it For any positive integer $m$ of the form (3.103), we have the wild solution (3.105) of the Erd\"{o}s-Straus equation.
If $m$ is of the form (3.109), we have the wild solution (3.111) of the Erd\"{o}s-Straus equation. }\psp

\section{Wild Solutions with $k=6$}

In this section, we solve the Er\"{o}s-Straus equation for $n=24m+1$ with $k=6$ in (3.1).  We divide it into the cases of even multiplier and odd multiplier, respectively.

\subsection{ Case of even multiplier}

In this subcase, we first consider
\begin{eqnarray}\frac{4}{n}&=&\frac{1}{6m+6}+\frac{23}{(6m+6)(24m+1)}\nonumber\\&=&
\frac{1}{6m+6}+\frac{2\ell\times 23}{2\ell\times6(m+1)(24m+1)}\nonumber\\&=&\frac{1}{6m+6}+\frac{(46\ell-1)+1}{12\ell(m+1)(24m+1)},\end{eqnarray}where $\ell$ is a positive integer.
Suppose $\ell\equiv 1\;(\mbox{mod}\;3)$; that is, $\ell=3r+1$ for some $r\in\mbb N$. Then
\begin{equation}46\ell-1=138r+45=3(46r+15).\end{equation}
According to (2.16),  we impose
\begin{equation}m+1\equiv 0\quad(\mbox{mod}\;46r+15).\end{equation}
Equivalently,
\begin{equation}m\equiv 46r+14\quad(\mbox{mod}\;46r+15).\end{equation}
Thus
\begin{equation}m=46r+14+c(46r+15)\qquad\mbox{for some}\;\;c\in\mbb N.\end{equation}

Assume that $m$ is of the above form. Then
\begin{equation}m+1=46r+14+c(46r+15)+1=(c+1)(46r+15).\end{equation}
According to (4.1),
\begin{eqnarray}\frac{4}{24m+1}&=&\frac{1}{6m+6}+\frac{(46\ell-1)+1}{12\ell(m+1)(24m+1)}\nonumber\\&=&
\frac{1}{6m+6}+\frac{3(46r+15)+1}{12(3r+1)(c+1)(46r+15)(24m+1)}\nonumber\\&=&
\frac{1}{6m+6}+\frac{1}{4(3r+1)(c+1)(24m+1)}\nonumber\\&&+\frac{1}{12(3r+1)(c+1)(46r+15)(24m+1)}.\end{eqnarray}
A subtle point is that we have well used the factor ``3" in $6m+6$ of the denominator. \psp

Next we assume $\ell\equiv 0\;(\mbox{mod}\;3)$; that is, $\ell=3r$ for some $0<r\in\mbb N$. Then
\begin{equation}46\ell-1=138r-1.\end{equation}
We assume
\begin{equation}m+1\equiv 0\quad(\mbox{mod}\;138r-1).\end{equation}
Equivalently,
\begin{equation}m\equiv 138r-2\quad(\mbox{mod}\;138r-1).\end{equation}
Thus
\begin{equation}m=138r-2+c(138r-1)\qquad\mbox{for some}\;\;c\in\mbb N.\end{equation}

When $m$ is of the above form,
\begin{equation}m+1=138r-2+c(138r-1)+1=(c+1)(138r-1).\end{equation}
By (4.1),
\begin{eqnarray}\frac{4}{24m+1}&=&\frac{1}{6m+6}+\frac{(46\ell-1)+1}{12\ell(m+1)(24m+1)}\nonumber\\&=&
\frac{1}{6m+6}+\frac{(138r-1)+1}{36r(c+1)(138r-1)(24m+1)}\nonumber\\&=&
\frac{1}{6m+6}+\frac{1}{36r(c+1)(24m+1)}\nonumber\\&&+\frac{1}{36r(c+1)(138r-1)(24m+1)}
.\end{eqnarray}
\pse

Thirdly we assume $\ell\equiv 2\;(\mbox{mod}\;3)$; that is, $\ell=3r+2$ for some $r\in\mbb N$. Then
\begin{equation}46\ell-1=138r+91.\end{equation}
We impose
\begin{equation}m+1\equiv 0\quad(\mbox{mod}\;138r+91).\end{equation}
Equivalently,
\begin{equation}m\equiv 138r+90\quad(\mbox{mod}\;138r+91).\end{equation}
Thus
\begin{equation}m=138r+90+c(138r+91)\qquad\mbox{for some}\;\;c\in\mbb N.\end{equation}

If $m$ is of the above form, then
\begin{equation}m+1=138r+90+c(138r+91)+1=(c+1)(138r+91).\end{equation}
Expression (4.1) yields
\begin{eqnarray}\frac{4}{24m+1}&=&\frac{1}{6m+6}+\frac{(46\ell-1)+1}{12\ell(m+1)(24m+1)}\nonumber\\&=&
\frac{1}{6m+6}+\frac{(138r+91)+1}{12(3r+2)(c+1)(138r+91)
(24m+1)}\nonumber\\&=&\frac{1}{6m+6}+\frac{1}{12(3r+2)(c+1)
(24m+1)}\nonumber\\&&+\frac{1}{12(3r+2)(c+1)(138r+91)
(24m+1)}.\end{eqnarray}
\pse

Now we  consider
\begin{eqnarray}\frac{4}{n}&=&\frac{1}{6m+6}+\frac{23}{(6m+6)(24m+1)}\nonumber\\&=&
\frac{1}{6m+6}+\frac{2\ell\times 23}{2\ell\times6(m+1)(24m+1)}\nonumber\\&=&\frac{1}{6m+6}+\frac{(46\ell-3)+3}{12\ell(m+1)(24m+1)},\end{eqnarray}where $\ell$ is a positive integer.
According to (1.8) with $\Im'_1=46\ell-3$ and $\Im'_2=3$, $\ell\not\equiv 0\;(\mbox{mod}\;3)$. Suppose $\ell\equiv 1\;(\mbox{mod}\;3)$; that is, $\ell=3r+1$ for some $r\in\mbb N$. Then
\begin{equation}46\ell-3=138r+43.\end{equation}
We assume
\begin{equation}m+1\equiv 0\quad(\mbox{mod}\;138r+43).\end{equation}
Equivalently,
\begin{equation}m\equiv 138r+42\quad(\mbox{mod}\;138r+43).\end{equation}
Thus
\begin{equation}m=138r+42+c(138r+43)\qquad\mbox{for some}\;\;c\in\mbb N.\end{equation}

When $m$ is of the above form,
\begin{equation}m+1=138r+42+c(138r+43)+1=(c+1)(138r+43).\end{equation}
Expression (4.20) gives
\begin{eqnarray}\frac{4}{24m+1}&=&\frac{1}{6m+6}+\frac{(46\ell-3)+3}{12\ell(m+1)(24m+1)}\nonumber\\&=&
\frac{1}{6m+6}+\frac{(138r+43)+3}{12(3r+1)(c+1)(138r+43)
(24m+1)}\nonumber\\&=&\frac{1}{6m+6}+\frac{1}{12(3r+1)(c+1)
(24m+1)}\nonumber\\&&+\frac{1}{4(3r+1)(c+1)(138r+43)
(24m+1)}.\end{eqnarray}
\pse

We assume $\ell\equiv 2\;(\mbox{mod}\;3)$; that is, $\ell=3r+2$ for some $r\in\mbb N$. Then
\begin{equation}46\ell-3=138r+89.\end{equation}
We impose
\begin{equation}m+1\equiv 0\quad(\mbox{mod}\;138r+89).\end{equation}
Equivalently,
\begin{equation}m\equiv 138r+88\quad(\mbox{mod}\;138r+89).\end{equation}
Thus
\begin{equation}m=138r+89+c(138r+89)\qquad\mbox{for some}\;\;c\in\mbb N.\end{equation}

If $m$ is of the above form, then
\begin{equation}m+1=138r+88+c(138r+89)+1=(c+1)(138r+89).\end{equation}
Expression (4.20) yields
\begin{eqnarray}\frac{4}{24m+1}&=&\frac{1}{6m+6}+\frac{(46\ell-3)+3}{12\ell(m+1)(24m+1)}\nonumber\\&=&
\frac{1}{6m+6}+\frac{(138r+89)+3}{12(3r+2)(c+1)(138r+89)
(24m+1)}\nonumber\\&=&\frac{1}{6m+6}+\frac{1}{12(3r+2)(c+1)
(24m+1)}\nonumber\\&&+\frac{1}{4(3r+2)(c+1)(138r+89)
(24m+1)}.\end{eqnarray}\pse

{\bf Theorem 4.1}\quad {\it For any positive integer $m$ of the form (4.5), we have the wild solution (4.7) of the Erd\"{o}s-Straus equation.
If $m$ is of the form (4.11), we have the wild solution (4.13) of the Erd\"{o}s-Straus equation. When $m$ is of the form (4.17), we have the wild solution (4.19) of the Erd\"{o}s-Straus equation. Suppose that $m$ is of the form (4.24), and then  we have the wild solution (4.26) of the Erd\"{o}s-Straus equation. Let $m$ be of the form (4.30), and then we have the wild solution (4.32) of the Erd\"{o}s-Straus equation.
}\psp

\subsection{ Case of odd multiplier}

In this case,  we firs consider
\begin{eqnarray}\frac{4}{n}&=&\frac{1}{6m+6}+\frac{23}{(6m+6)(24m+1)}\nonumber\\&=&
\frac{1}{6m+6}+\frac{(2\ell+1)\times 23}{(2\ell+1)\times6(m+1)(24m+1)}\nonumber\\&=&\frac{1}{6m+6}+\frac{2(23\ell+11)+1}{6(2\ell+1)(m+1)(24m+1)},\end{eqnarray}where $\ell$ is a positive integer.
Suppose $\ell\equiv 2\;(\mbox{mod}\;3)$; that is, $\ell=3r+2$ for some $r\in\mbb N$. Then
\begin{equation}23\ell+11=69r+57=3(23r+19).\end{equation}
By (2.16),  we impose
\begin{equation}m+1\equiv 0\quad(\mbox{mod}\;23r+19).\end{equation}
Equivalently,
\begin{equation}m\equiv 23r+18\quad(\mbox{mod}\;23r+19).\end{equation}
Thus
\begin{equation}m=23r+18+c(23r+19)\qquad\mbox{for some}\;\;c\in\mbb N.\end{equation}

Assume that $m$ is of the above form. Then
\begin{equation}m+1=23r+18+c(23r+19)+1=(c+1)(23r+19).\end{equation}
According to (4.33),
\begin{eqnarray}\frac{4}{24m+1}&=&\frac{1}{6m+6}+\frac{2(23\ell+11)+1}{6(2\ell+1)(m+1)(24m+1)}
\nonumber\\&=&\frac{1}{6m+6}+\frac{6(23r+19)+1}{6(6r+5)(c+1)(23r+19)(24m+1)}
\nonumber\\&=&\frac{1}{6m+6}+\frac{1}{(6r+5)(c+1)(24m+1)}
\nonumber\\&&+\frac{1}{6(6r+5)(c+1)(23r+19)(24m+1)}.\end{eqnarray}
A subtle point is that we have well used the factor ``6" in $6m+6$ of the denominator. \psp

Next we assume $\ell\equiv 0\;(\mbox{mod}\;3)$; that is, $\ell=3r$ for some $0<r\in\mbb N$. Then
\begin{equation}23\ell+11=69r+11.\end{equation}
We assume
\begin{equation}m+1\equiv 0\quad(\mbox{mod}\;69r+11).\end{equation}
Equivalently,
\begin{equation}m\equiv 69r+10\quad(\mbox{mod}\;69r+11).\end{equation}
Thus
\begin{equation}m=69r+10+c(69r+11)\qquad\mbox{for some}\;\;c\in\mbb N.\end{equation}

When $m$ is of the above form,
\begin{equation}m+1=69r+10+c(69r+11)+1=(c+1)(69r+11).\end{equation}
By (4.33),
\begin{eqnarray}\frac{4}{24m+1}&=&\frac{1}{6m+6}+\frac{2(23\ell+11)+1}{6(2\ell+1)(m+1)(24m+1)}
\nonumber\\&=&\frac{1}{6m+6}+\frac{2(69r+11)+1}{6(6r+1)(c+1)(69r+11)(24m+1)}
\nonumber\\&=&\frac{1}{6m+6}+\frac{1}{3(6r+1)(c+1)(24m+1)}\nonumber\\&&+\frac{1}{6(6r+1)(c+1)(69r+11)(24m+1)}.\end{eqnarray}
\pse

Let $\ell\equiv 1\;(\mbox{mod}\;3)$; that is, $\ell=3r+1$ for some $0<r\in\mbb N$. Then
\begin{equation}23\ell+11=69r+34.\end{equation}
We assume
\begin{equation}m+1\equiv 0\quad(\mbox{mod}\;69r+34).\end{equation}
Equivalently,
\begin{equation}m\equiv 69r+33\quad(\mbox{mod}\;69r+34).\end{equation}
Thus
\begin{equation}m=69r+33+c(69r+34)\qquad\mbox{for some}\;\;c\in\mbb N.\end{equation}

If $m$ be of the above form, then
\begin{equation}m+1=69r+33+c(69r+34)+1=(c+1)(69r+34).\end{equation}
By (4.43),
\begin{eqnarray}\frac{4}{24m+1}&=&\frac{1}{6m+6}+\frac{2(23\ell+11)+1}{6(2\ell+1)(m+1)(24m+1)}
\nonumber\\&=&\frac{1}{6m+6}+\frac{2(69r+34)+1}{18(2r+1)(c+1)(69r+34)(24m+1)}
\nonumber\\&=&\frac{1}{6m+6}+\frac{1}{9(2r+1)(c+1)(24m+1)}
\nonumber\\&&+\frac{1}{18(2r+1)(c+1)(69r+34)(24m+1)}.\end{eqnarray}
\pse

Furthermore, we consider
\begin{eqnarray}\frac{4}{n}&=&\frac{1}{6m+6}+\frac{23}{(6m+6)(24m+1)}\nonumber\\&=&
\frac{1}{6m+6}+\frac{(2\ell+1)\times 23}{(2\ell+1)\times6(m+1)(24m+1)}\nonumber\\&=&\frac{1}{6m+6}+\frac{(46\ell+21)+2}{6(2\ell+1)(m+1)(24m+1)},\end{eqnarray}where $\ell$ is a positive integer.
Suppose $\ell\equiv 0\;(\mbox{mod}\;3)$; that is, $\ell=3r$ for some $0<r\in\mbb N$. Then
\begin{equation}46\ell+21=3(46r+7).\end{equation}
By (2.16),  we impose
\begin{equation}m+1\equiv 0\quad(\mbox{mod}\;46r+7).\end{equation}
Equivalently,
\begin{equation}m\equiv 46r+6\quad(\mbox{mod}\;46r+7).\end{equation}
Thus
\begin{equation}m=46r+6+c(46r+7)\qquad\mbox{for some}\;\;c\in\mbb N.\end{equation}

Assume that $m$ is of the above form. Then
\begin{equation}m+1=46r+7+c(46r+7)+1=(c+1)(46r+7).\end{equation}
According to (4.52),
\begin{eqnarray}\frac{4}{24m+1}&=&\frac{1}{6m+6}+\frac{(46\ell+21)+2}{6(2\ell+1)(m+1)(24m+1)}
\nonumber\\&=&\frac{1}{6m+6}+\frac{3(46r+7)+2}{6(6r+1)(c+1)(46r+7)(24m+1)}
\nonumber\\&=&\frac{1}{6m+6}+\frac{1}{2(6r+1)(c+1)(24m+1)}
\nonumber\\&&+\frac{1}{3(6r+1)(c+1)(46r+7)(24m+1)}.\end{eqnarray}
A subtle point is that we have well used the factor ``3" in $6m+6$ of the denominator. \psp

Next we assume $\ell\equiv 1\;(\mbox{mod}\;3)$; that is, $\ell=3r+1$ for some $r\in\mbb N$. Then
\begin{equation}46\ell+21=138r+67.\end{equation}
So we impose
\begin{equation}m+1\equiv 0\quad(\mbox{mod}\;138r+67).\end{equation}
Equivalently,
\begin{equation}m\equiv 138r+66\quad(\mbox{mod}\;138r+67).\end{equation}
Thus
\begin{equation}m=138r+66+c(138r+67)\qquad\mbox{for some}\;\;c\in\mbb N.\end{equation}

Let $m$ be of the above form. Then
\begin{equation}m+1=138r+66+c(138r+67)+1=(c+1)(138r+67).\end{equation}
By (4.52),
\begin{eqnarray}\frac{4}{24m+1}&=&\frac{1}{6m+6}+\frac{(46\ell+21)+2}{6(2\ell+1)(m+1)(24m+1)}
\nonumber\\&=&\frac{1}{6m+6}+\frac{(138r+67)+2}{18(2r+1)(c+1)(138r+67)(24m+1)}
\nonumber\\&=&\frac{1}{6m+6}+\frac{1}{18(2r+1)(c+1)(24m+1)}
\nonumber\\&&+\frac{1}{9(2r+1)(c+1)(138r+67)(24m+1)}.\end{eqnarray}
When $r=0$ and $c=442$, the above equation becomes (2.18).
\psp

Let $\ell\equiv 2\;(\mbox{mod}\;3)$; that is, $\ell=3r+2$ for some $r\in\mbb N$. Then
\begin{equation}46\ell+21=138r+113.\end{equation}
So we impose
\begin{equation}m+1\equiv 0\quad(\mbox{mod}\;138r+113).\end{equation}
Equivalently,
\begin{equation}m\equiv 138r+112\quad(\mbox{mod}\;138r+113).\end{equation}
Hence
\begin{equation}m=138r+112+c(138r+113)\qquad\mbox{for some}\;\;c\in\mbb N.\end{equation}

Assume that $m$ is of the above form. Then
\begin{equation}m+1=138r+112+c(138r+113)+1=(c+1)(138r+113).\end{equation}
Expression (4.52) yields
\begin{eqnarray}\frac{4}{24m+1}&=&\frac{1}{6m+6}+\frac{(46\ell+21)+2}{6(2\ell+1)(m+1)(24m+1)}
\nonumber\\&=&\frac{1}{6m+6}+\frac{(138r+113)+2}{6(6r+5)(c+1)(138r+113)(24m+1)}
\nonumber\\&=&\frac{1}{6m+6}+\frac{1}{6(6r+5)(c+1)(24m+1)}
\nonumber\\&&+\frac{1}{3(6r+5)(c+1)(138r+113)(24m+1)}.\end{eqnarray}\pse

In addition, we consider
\begin{eqnarray}\frac{4}{n}&=&\frac{1}{6m+6}+\frac{23}{(6m+6)(24m+1)}\nonumber\\&=&
\frac{1}{6m+6}+\frac{(2\ell+1)\times 23}{(2\ell+1)\times6(m+1)(24m+1)}\nonumber\\&=&\frac{1}{6m+6}+\frac{2(23\ell+10)+3}{6(2\ell+1)(m+1)(24m+1)},\end{eqnarray}where $\ell$ is a positive integer. According to (1.8) with $\Im'_1=2(23\ell+10)+3$ and $\Im'_2=3$, $\ell\not\equiv 1\;(\mbox{mod}\;3)$.
Suppose $\ell\equiv 0\;(\mbox{mod}\;3)$; that is, $\ell=3r$ for some $0<r\in\mbb N$. Then
\begin{equation}23\ell+10=69r+10.\end{equation}
So we impose
\begin{equation}m+1\equiv 0\quad(\mbox{mod}\;69r+10).\end{equation}
Equivalently,
\begin{equation}m\equiv 69r+9\quad(\mbox{mod}\;69r+10).\end{equation}
Thus
\begin{equation}m=69r+9+c(69r+10)\qquad\mbox{for some}\;\;c\in\mbb N.\end{equation}

Assume that $m$ is of the above form. Then
\begin{equation}m+1=69r+9+c(69r+10)+1=(c+1)(69r+10).\end{equation}
According to (4.71),
\begin{eqnarray}\frac{4}{24m+1}&=&\frac{1}{6m+6}+\frac{2(23\ell+10)+3}{6(2\ell+1)(m+1)(24m+1)}\nonumber\\&=&
\frac{1}{6m+6}+\frac{2(69r+10)+3}{6(6r+1)(c+1)(69r+10)(24m+1)}\nonumber\\&=&
\frac{1}{6m+6}+\frac{1}{3(6r+1)(c+1)(24m+1)}\nonumber\\&&+\frac{1}{2(6r+1)(c+1)(69r+10)(24m+1)}
.\end{eqnarray}
A subtle point is that we have well used the factor ``3" in $6m+6$ of the denominator. \psp

Let $\ell\equiv 2\;(\mbox{mod}\;3)$; that is, $\ell=3r+2$ for some $0<r\in\mbb N$. Then
\begin{equation}23\ell+10=69r+56.\end{equation}
So we impose
\begin{equation}m+1\equiv 0\quad(\mbox{mod}\;69r+56).\end{equation}
Equivalently,
\begin{equation}m\equiv 69r+55\quad(\mbox{mod}\;69r+56).\end{equation}
Thus
\begin{equation}m=69r+55+c(69r+56)\qquad\mbox{for some}\;\;c\in\mbb N.\end{equation}

When $m$ is of the above form,
\begin{equation}m+1=69r+55+c(69r+56)+1=(c+1)(69r+56).\end{equation}
By (4.71),
\begin{eqnarray}\frac{4}{24m+1}&=&\frac{1}{6m+6}+\frac{2(23\ell+10)+3}{6(2\ell+1)(m+1)(24m+1)}\nonumber\\&=&
\frac{1}{6m+6}+\frac{2(69r+56)+3}{6(6r+5)(c+1)(69r+56)(24m+1)}\nonumber\\&=&
\frac{1}{6m+6}+\frac{1}{3(6r+5)(c+1)(24m+1)}\nonumber\\&&+\frac{1}{2(6r+5)(c+1)(69r+56)(24m+1)}
.\end{eqnarray}
\pse

Next we we assume that $m=2m_1+1$ with $m_1\in\mbb N$. We consider
\begin{eqnarray}\frac{4}{n}&=&\frac{1}{6m+6}+\frac{23}{(6m+6)(24m+1)}\nonumber\\&=&
\frac{1}{6m+6}+\frac{(2\ell+1)\times 23}{(2\ell+1)\times6(m+1)(24m+1)}\nonumber\\&=&\frac{1}{6m+6}+\frac{(46\ell+19)+4}{12(2\ell+1)(m_1+1)(24m+1)},\end{eqnarray}where $\ell$ is a positive integer. First we suppose
\begin{equation}(46\ell+19)|[3(m_1+1)].\end{equation}
Assume $\ell=3s+2$. Then
\begin{equation}
46\ell+19=138s+111=3(46s+37)\lra (46s+37)|(m_1+1).\end{equation}
So
\begin{equation}m_1=t(46s+37)+46s+36=36+37t+46s(t+1).\end{equation}
In particular,
\begin{equation}m=73+74t+92s(t+1).\end{equation}
By (4.84),
\begin{eqnarray}\frac{4}{n}&=&\frac{1}{6m+6}+\frac{(46\ell+19)+4}{12(2\ell+1)(m_1+1)(24m+1)}\nonumber\\&=&
\frac{1}{6m+6}+\frac{3(46s+37)+4}{12(6s++5)(t+1)(46s+37)(24m+1)}\nonumber\\&=&
\frac{1}{6m+6}+\frac{1}{4(6s++5)(t+1)(24m+1)}\nonumber\\&&+\frac{1}{3(6s++5)(t+1)(46s+37)(24m+1)}.\end{eqnarray}
When $m=71925$, $n=24m+1=1726201$ is a prime. Taking $s=781$ and $t=0$, we have $2\ell+1=4691$ and
\begin{equation}\frac{1}{1726201}=\frac{1}{431556}+\frac{1}{32390435564}+\frac{1}{873642925641099},\end{equation}
which coincides with (2.23).

Assume $3\not|(46\ell+19)$. Then
\begin{equation}(46\ell+19)|(m_1+1)\lra m_1=t(46\ell+19)+46\ell+18=18+19t+46\ell(t+1)
.\end{equation}
In particular,
\begin{equation}m=37+38t+92\ell(t+1)).\end{equation}
By (4.84),
\begin{eqnarray}\frac{4}{n}&=&\frac{1}{6m+6}+\frac{(46\ell+19)+4}{12(2\ell+1)(m_1+1)(24m+1)}\nonumber\\&=&
\frac{1}{6m+6}+\frac{(46\ell+19)+4}{12(2\ell+1)(t+1)(46\ell+19)(24m+1)}\nonumber\\&=&
\frac{1}{6m+6}+\frac{1}{12(2\ell+1)(t+1)(24m+1)}\nonumber\\&&+\frac{1}{3(2\ell+1)(t+1)(46\ell+19)(24m+1)}.\end{eqnarray}
\pse

Finally we consider
\begin{eqnarray}\frac{4}{n}&=&\frac{1}{6m+6}+\frac{23}{(6m+6)(24m+1)}\nonumber\\&=&
\frac{1}{6m+6}+\frac{(2\ell+1)\times 23}{(2\ell+1)\times6(m+1)(24m+1)}\nonumber\\&=&\frac{1}{6m+6}+\frac{(46\ell+17)+6}{6(2\ell+1)(m+1)(24m+1)},\end{eqnarray}where $\ell$ is a positive integer. According to (1.8) with $\Im'_1=46\ell+17$ and $\Im'_2=3$, $\ell\not\equiv 1\;(\mbox{mod}\;3)$.
Suppose $\ell\equiv 0\;(\mbox{mod}\;3)$; that is, $\ell=3r$ for some $0<r\in\mbb N$. Then
\begin{equation}46\ell+17=138r+17.\end{equation}
By (2.13), we impose
\begin{equation}m+1\equiv 0\quad(\mbox{mod}\;138r+17).\end{equation}
Equivalently,
\begin{equation}m\equiv 138r+16\quad(\mbox{mod}\;138r+17).\end{equation}
Thus
\begin{equation}m=138r+16+c(138r+17)\qquad\mbox{for some}\;\;c\in\mbb N.\end{equation}

Assume that $m$ is of the above form. Then
\begin{equation}m+1=138r+16+c(138r+17)+1=(c+1)(138r+17).\end{equation}
According to (4.84),
\begin{eqnarray}\frac{4}{24m+1}&=&\frac{1}{6m+6}+\frac{(46\ell+17)+6}{6(2\ell+1)(m+1)(24m+1)}
\nonumber\\&=&\frac{1}{6m+6}+\frac{(138r+17)+6}{6(6r+1)(c+1)(138r+17)(24m+1)}
\nonumber\\&=&\frac{1}{6m+6}+\frac{1}{6(6r+1)(c+1)(24m+1)}
\nonumber\\&&+\frac{1}{(6r+1)(c+1)(138r+17)(24m+1)}.\end{eqnarray}
\pse

Let $\ell\equiv 2\;(\mbox{mod}\;3)$; that is, $\ell=3r+2$ for some $0<r\in\mbb N$. Then
\begin{equation}46\ell+17=138r+109.\end{equation}
By (2.13),  we impose
\begin{equation}m+1\equiv 0\quad(\mbox{mod}\;138r+109).\end{equation}
Equivalently,
\begin{equation}m\equiv 138r+108\quad(\mbox{mod}\;138r+109).\end{equation}
Thus
\begin{equation}m=138r+108+c(138r+109)\qquad\mbox{for some}\;\;c\in\mbb N.\end{equation}

Assume that $m$ is of the above form. Then
\begin{equation}m+1=138r+108+c(138r+109)+1=(c+1)(138r+109).\end{equation}
By (4.84),
\begin{eqnarray}\frac{4}{24m+1}&=&\frac{1}{6m+6}+\frac{(46\ell+17)+6}{6(2\ell+1)(m+1)(24m+1)}
\nonumber\\&=&\frac{1}{6m+6}+\frac{(138r+109)+6}{6(6r+5)(c+1)(138r+109)(24m+1)}
\nonumber\\&=&\frac{1}{6m+6}+\frac{1}{6(6r+5)(c+1)(24m+1)}
\nonumber\\&&+\frac{1}{(6r+5)(c+1)(138r+109)(24m+1)}.\end{eqnarray}
\pse

{\bf Theorem 4.2}\quad {\it For any positive integer $m$ of the form (4.37), we have the wild solution (4.39) of the Erd\"{o}s-Straus equation.
If $m$ is of the form (4.43), we have the wild solution (4.45) of the Erd\"{o}s-Straus equation. When $m$ is of the form (4.49), we have the wild solution (4.51) of the Erd\"{o}s-Straus equation. Suppose that $m$ is of the form (4.56), and then  we have the wild solution (4.58) of the Erd\"{o}s-Straus equation. Let $m$ be of the form (4.62), and then we have the wild solution (4.64) of the Erd\"{o}s-Straus equation.

For any positive integer $m$ of the form (4.88), we have the wild solution (4.89) of the Erd\"{o}s-Straus equation. For any positive integer $m$ of the form (4.92), we have the wild solution (4.93) of the Erd\"{o}s-Straus equation.
For any positive integer $m$ of the form (4.68), we have the wild solution (4.70) of the Erd\"{o}s-Straus equation.
If $m$ is of the form (4.75), we have the wild solution (4.77) of the Erd\"{o}s-Straus equation. When $m$ is of the form (4.81), we have the wild solution (4.83) of the Erd\"{o}s-Straus equation. Suppose that $m$ is of the form (4.98), and then we have the wild solution (4.100) of the Erd\"{o}s-Straus equation. Let $m$ be of the form (4.104), and then we have the wild solution (4.1066) of the Erd\"{o}s-Straus equation.}\psp

\pse

\psp

E-Mail:  xiaoping@math.ac.cn

 \end{document}